\documentclass[a4paper,12pt]{article}
\usepackage[cp850]{inputenc} 
\title{Pyramids and \mbus}
\author{M.J. Soto\\ \tt soto@us.es \\ Departamento de Algebra\\
Universidad de Sevilla \and Jos\'{e} L. Vicente\\ \tt jlvc@us.es\\
Departamento de Algebra\\ Universidad de Sevilla
}
\usepackage{amsmath}
\usepackage{amssymb}
\makeatletter
\usepackage[dvips]{graphicx}
\DeclareMathAlphabet{\bboard}{U}{msb}{m}{n}

\def\RR{\bboard{R}}
\def\ZZ{\bboard{Z}}

  \newcounter{numero}[section]
  \newcounter{aux1}
  \newcounter{aux2}
  \renewcommand{\thenumero}{\thesection.\arabic{numero}}
  \newenvironment{corolario}{
  \dimen255=\parindent \parindent=0in
  \refstepcounter{numero}{\vspace{.3cm}\bf Corollary \thenumero.--}\
  \it}{\rm \parindent=\dimen255\par}

  \newenvironment{definicion}{
  \dimen255=\parindent \parindent=0in
  \refstepcounter{numero}{\vspace{.3cm}\bf Definition
\thenumero.--}\
  \it}{\rm \parindent=\dimen255\par}

  \newenvironment{teorema}{
  \dimen255=\parindent \parindent=0in
  \refstepcounter{numero}{\vspace{.3cm}\bf Theorem
\thenumero.--}\
  \it}{\rm \parindent=\dimen255\par}

  \newenvironment{lema}{
  \dimen255=\parindent \parindent=0in
  \refstepcounter{numero} {\vspace{.3cm}\bf Lemma
\thenumero.--}\
  \it}{\rm \parindent=\dimen255\par}

  \newenvironment{nota}{
  \refstepcounter{numero} {\vspace{.3cm}\noindent\bf Remark
\thenumero.--}\
  }{
  \par}

  \newenvironment{proposicion}{
  \dimen255=\parindent \parindent=0in
  \refstepcounter{numero} {\vspace{.3cm}\bf Proposition
\thenumero.--}\ \it
  }{\rm \parindent=\dimen255\par}

  \newenvironment{notas}{\setcounter{aux1}{0}
  \dimen255=\parindent \parindent=0in
  \refstepcounter{numero} {\vspace{.3cm}\bf Remarks
\thenumero.--}\
  }{\parindent=\dimen255\par}

  \newcommand{\parnot}{\vspace{.2cm} \parindent=0in
  \refstepcounter{aux1}
  \thenumero.\arabic{aux1}.\hspace{.1cm}\parindent=\dimen255}

\def\pushright#1{{
   \parfillskip=0pt            
   \widowpenalty=10000         
   \displaywidowpenalty=10000  
   \finalhyphendemerits=0      
   \leavevmode                 
   \unskip                     
   \nobreak                    
   \hfil                       
   \penalty50                  
   \hskip.2em                  
   \null                       
   \hfill                      
   {#1}                        
   \par}}                      

\def\qed{\pushright{\framebox[2.5mm]{}}\penalty-700 \smallskip}

  \newenvironment{demostracion}{
  {\vspace{.3cm}\noindent\bf{Proof:}}\
  }{\qed\par}

\def\[[#1]]{[\![#1]\!]}

\def\pe#1{\langle #1 \rangle}
\let\em=\it

\def\mlex{<_{\hbox{\scriptsize lex}}}

\def\ai{affi\-ne\-ly in\-de\-pen\-dent}

\def\vg#1#2{(#1_1,\ldots,#1_{#2})}

\def\dps{\colon\!}
\def\nd{New\-ton dia\-gram}
\def\nds{New\-ton dia\-grams}
\def\mbu{mo\-no\-mial blo\-wing-up}
\def\mbus{mo\-no\-mial blo\-wing-ups}
\def\mbd{mo\-no\-mial blo\-wing\--down}
\def\mbds{mo\-no\-mial blo\-wing\--downs}
\def\ncd{nor\-mal cros\-sing di\-vi\-sor}

\def\pc{first qua\-drant}
\def\opc{op\-po\-si\-te of the first qua\-drant}
\def\se{half\--spa\-ce}

\def\tsr{half\--li\-ne}
\def\tsrs{half\--li\-nes}
\def\sr#1{\pe{#1}_{+}}

\newcommand{\tmcd}{gre\-at\-est com\-mon di\-vi\-sor}

\makeatother
\begin{document}

\maketitle

\abstract{We show that a convex pyramid $\Gamma$ in
$\RR^{n}$ with
apex at $\mathbf{0}$ can be brought to the \pc\ by a finite
sequence of \mbus\  if and only if
$\Gamma\cap(-\RR_{0}^{n})=\{\mathbf{0}\}$. The proof is
non-trivially derived from the theorem of Farkas-Minkowski.
Then, we
apply this theorem to show how the \nds\ of the roots of any
Weierstra\ss\ polynomial
$$
P(\mathbf{x},z)=z^{m}+h_{1}(\mathbf{x})z^{m-1}+\cdots+h_{m-1}(\mathbf{x})z+h_{m}(\mathbf{x})\, ,
$$
$h_{i}(\mathbf{x})\in k\[[x_{1},\ldots,x_{n}]][z]$,
are contained in a pyramid of this type. Finally, if $n=2$,
this fact is equivalent to the Jung-Abhyankar theorem.}

\section{Introduction} We will operate in the
euclidean space $\RR^{n}$ with its affine structure. As it
is classical, we
will distinguish between the point-space $X=\RR^{n}$ and the
underlying vector space $V=\RR^{n}$. The vector addition
is the canonical action (translations) of $V$ on $X$.

A {\em polyhedron}\index{polyhedron} is the
convex hull of a finite set of points generating the affine
space
$X$. Equivalently (\cite{Ewald}, page 30),
a polyhedron is
the compact intersection of a finite
set of half-spaces. Moreover, if
$E=\{A_{0},A_{1},\ldots,A_{m}\}$ is a finite set of points
generating $X$  and, if $\mathcal{H}=\{H_{1},\ldots,H_{p}\}$
is the set
of all (different) hyperplanes passing through all possible
subsets of $E$ consisting of $n$ affinely independent
points, then the set of vertices of the convex hull $[E]$ of
$E$ is the set of points of intersection of all the subsets
of $\mathcal{H}$ consisting of $n$ hyperplanes whose
intersection is an only point\footnote{c.f., for instance, Vicente, J.L.:
{\em Notas sobre convexidad} at http://www.us.es/da}.

A {\em pyramid}\index{pyramid} is the
projection, from the origin $\mathbf{0}\in\RR^{n}$, of a
polyhedron contained in a hyperplane $H$ not passing through
$\mathbf{0}$. To be more precise: given a hyperplane $H$
such that
$\mathbf{0}\notin H$ and a finite set $E$ of points
generating
$H$,
denoting by $\Delta$ the convex hull $[E]$ of $E$, then the
corresponding pyramid is
$$
\Gamma(\Delta)=\bigcup_{\mathbf{a}\in\Delta}\pe{\mathbf{a}}_{+}\,
,
$$
where $\sr{\mathbf{a}}$ is the \tsr\ of the non-negative
multiples of $\mathbf{a}$.
Equivalently (c.f. Vicente, J.L., loc.cit.),
a pyramid can be given by a polyhedron $\Delta'$
in $X$, one of whose vertices is $\mathbf{0}$. In this case,
the pyramid is nothing but
$$
\Gamma=\bigcup_{\mathbf{a}\in\Delta'\setminus\{O\}}
\sr{\mathbf{a}}\,
.
$$
Moreover, if
$\{\mathbf{0},\mathbf{a}_{1},\ldots,\mathbf{a}_{m}\}$
are the vertices of
$\Delta'$, there is a hyperplane $H$ strictly separating
$\mathbf{0}$
from $[\mathbf{a}_{1},\ldots,\mathbf{a}_{m}]$. Then
$\Delta=H\cap\Delta'$
is a polyhedron in $H$ and $\Gamma=\Gamma(\Delta)$.

\begin{definicion}{\label{0207041}}
We will call a
\mbu\index{\mbu} (resp. a \mbd\index{mbd}) any
linear automorphism of
$\RR^{n}$ of the form
$$
\vg{a}{n}\to\vg{a}{n}M_{ij}\quad(\hbox{resp.}\quad\vg{a}{n}\to\vg{a}{n}N_{ij})
$$
where:
\begin{enumerate}

\item $i,j\in\ZZ$, $i\neq j$, $1\leq i,j\leq n$

\item $M_{ij}$ (resp. $N_{ij}$) is equal to the identity
matrix in which the
$(i,j)$-entry is set to $1$ (resp. to $-1$).

\end{enumerate}
\end{definicion}

\begin{nota}{\label{0307041}} With the notations of
definition \ref{0207041}, the \mbu\ (resp. \mbd\ ) acts in
the following way:
\begin{eqnarray*}
&&
\vg{a}{n}
\to
(a_{1},\ldots,a_{i}\stackrel{j)}{+}a_{j},\ldots,a_{n})\\
\hbox{resp.}
&&
\vg{a}{n}
\to
(a_{1},\ldots,a_{j}\stackrel{j)}{-}a_{i},\ldots,a_{n})
\end{eqnarray*}
This corresponds to the behavior of the exponents of a
monomial under the geometric \mbu\ $x_{i}\to x_{i}x_{j}$
or the \mbd\ $x_{i}\to x_{i}/x_{j}$. In fact, this
geometric
\mbu\ (resp. \mbd\ )  acts on monomials in the following
way:
\begin{eqnarray*}
&&
x_{1}^{a_{1}}\cdots x_{j}^{a_{j}}\cdots x_{n}^{a_{n}}
\to
x_{1}^{a_{1}}\cdots x_{j}^{a_{i}+a_{j}}\cdots
x_{n}^{a_{n}} \\
\hbox{resp.}
&&
x_{1}^{a_{1}}\cdots x_{j}^{a_{j}}\cdots
x_{n}^{a_{n}}
\to
x_{1}^{a_{1}}\cdots x_{j}^{a_{j}-a_{i}}\cdots
x_{n}^{a_{n}} \, .
\end{eqnarray*}
This is the reason of the name for these linear
automorphisms.

{\em From now on, we will use the name of \mbu\ (resp. \mbd)
indistinctly for the linear automorphisms defined in
\ref{0207041} or for the polynomial substitutions}
$x_{i}\to x_{i}x_{j}$
{\em (resp.}$x_{i}\to x_{i}/x_{j}${\em )}.
\end{nota}

\begin{nota}{\label{0307042}}
Let $E=\{\mathbf{a}_{1},\ldots,\mathbf{a}_{m}\}$
be a finite set of points generating a hyperplane $H$ which
does not contain the origin,  let $\Delta=[E]\subset H$ be
the
corresponding polyhedron and $\Gamma(\Delta)$ the pyramid;
then
$$
\Gamma(\Delta)=\left\{\left.\sum_{i=1}^{m}\lambda_{i}\mathbf{a}_{i}
\right| \lambda_{i}\geq 0\, , \forall i=1,\ldots,m\right\}\,
.
$$
Equivalently, let $A$ be the matrix whose row
vectors are
$\{\mathbf{a}_{1},\ldots,\mathbf{a}_{m}\}$; then
$$
\Gamma(\Delta)=\left\{(\lambda_{1},\ldots,\lambda_{m})A\mid
(\lambda_{1},\ldots,\lambda_{m})\in\RR_{0}^{m}
\right\}\, .
$$
If $\varphi_{ij}$ is the \mbu\ (resp. the \mbd) with matrix
$M_{ij}$ (resp. $N_{ij}$) and
$\mathbf{a}'_{i}=\varphi_{ij}(\mathbf{a}_{i})$ then
$E'=\{\mathbf{a}'_{1},\ldots,\mathbf{a}'_{m}\}$
generates a hyperplane $H'$ not passing through
$\mathbf{0}$. If
$\Delta'=[E']$, then
$\Gamma(\Delta')=\varphi(\Gamma(\Delta))$,
so it makes sense to
speak on the transform of a pyramid by a \mbu\ or a \mbd.
\end{nota}

\begin{definicion}{\label{0307043}}
The
\pc\index{\pc} of $X$ is the set $\RR_{0}^{n}$.
The
\opc\index{\opc} of $X$ is the set
$-\RR_{0}^{n}$.
\end{definicion}

\vspace{0.3cm}

The main problem we deal with in this paper is whether,
given a pyramid $\Gamma(\Delta)$, it exists a finite
sequence of \mbus\ such that the transform of the pyramid by
the sequence is contained  in the \pc.
We solve it by giving a geometrical criterion, from which we
derive explicit computations using existing optimization
algorithms.
The criterion is the following:

\begin{teorema}{\label{0307044}}
Let $E=\{\mathbf{a}_{1},\ldots,\mathbf{a}_{m}\}$
be a finite set of points generating a hyperplane $H$
not containing the origin,  let $\Delta=[E]$ be the
corresponding polyhedron and $\Gamma(\Delta)$ the pyramid;
then the following conditions are equivalent:
\begin{enumerate}

\item There exists a finite sequence of \mbus\ such that the
transform of $\Gamma(\Delta)$ by the sequence is contained
in the \pc.

\item
$\Gamma(\Delta)\cap(-\RR_{0}^{n})=\{\mathbf{0}\}$.

\end{enumerate}
\end{teorema}

\vspace{0.3cm}

The second condition can be easily checked by
the simplex method; in remark \ref{0507041} we will show how.
Let us observe that the first condition implies the second
because
$-\RR_{0}^{n}$
is stable by \mbus.
In fact, if $\Gamma(\Delta)$ had a
point $\mathbf{a}\neq \mathbf{0}$ in common with
$-\RR_{0}^{n}$ then, no matter what
sequence of \mbus\ we apply, the transform of $\mathbf{a}$
will stay
in $-\RR_{0}^{n}$.  The point is then to prove that the
second condition implies the first.

As an application, we show that the \nds\  of all the
roots of a Weierstra\ss\ polynomial
$$
P(\mathbf{x},z)=z^{m}+h_{1}(\mathbf{x})z^{m-1}+\cdots+h_{m-1}(\mathbf{x})z+h_{m}(\mathbf{x})
\in k\[[\mathbf{x}]][z] \, , \quad m>1
$$
are contained in a pyramid satifying the equivalent
conditions of theorem  \ref{0307044}. Moreover, we show how,
in dimension $2$, this fact is in some sense
equivalent to the Jung-Abhyankar theorem (c.f.
\cite{Abhyankar}).

\section{The proof}
\begin{nota}{\label{0407044}} Let us denote by $A$ the
matrix whose row vectors are
$\{\mathbf{a}_{1},\ldots,\mathbf{a}_{m}\}$. We will
speak of ``bringing $A$ to the \pc'' as equivalent to
bringing $\Gamma(\Delta)$ to the \pc. It is obvious that, if
$A$ has  column with only positive entries, then it can be
brought to the \pc by a finite sequence of \mbus: it is
enough
to add this column to the others a suitable number of times.
\end{nota}

\vspace{0.3cm}

We use linear optimization methods to prove theorem \ref{0307044}.
We refer to \cite{Gale} for the theorem
of Farkas-Minkowski and its consequences. In particular, we
take from it (pages 42-51) the following consequence of this
theorem (which might be also taken as an alternative
statement of it):

\begin{corolario}{\label{0407045}} Let $A$ be a matrix
$m\times n$; then one and only one of the following
conditions hold:
\begin{enumerate}

\item There exists a non-zero vector $x\geq 0$ such
that
$xA\leq 0$.

\item The system of inequalities $Ay>0$ has a non-negative
solution.

\end{enumerate}
Inequalities must be understood componentwise.
\end{corolario}

\vspace{0.3cm}

From this result we derive the following, which is the
useful one:

\begin{corolario}{\label{0307045}} Let $A$  be a matrix
$m\times n$; then one and only one of the following
conditions holds:
\begin{enumerate}

\item There exists a non-zero vector $x\geq 0$ such that
$xA\leq 0$.

\item The system of inequalities $Ay>0$ has a positive
integer solution (that is, a vector of positive integers).

\end{enumerate}
\end{corolario}

\begin{demostracion}
By corollary \ref{0407045}, we must only prove
that the existence of a non-negative real solution of
$Ay>0$ implies that there is a positive integer one.
Hence, it is enough to show the existence of a positive
rational solution.

Let $y$ be a non-negative real column vector such that
$Ay>0$, let
$a_i$ be the rows of $A$, $i=1,\ldots,m$, and let
$a_iy=\alpha_i>0$. Let $0<\varepsilon<\min\{\alpha_i\}$ and
let $K$ be a polidisc centered at $y$ such that, for all
$z\in K$, one has $|-\alpha_i+a_iz|<\varepsilon$. Then,
$-\varepsilon<a_iz-\alpha_i<\varepsilon$, so
$\alpha_i-\varepsilon<a_iz<\varepsilon+\alpha_i$,
which implies that $a_iz>0$. Since $K$ contains rational
vectors greater than zero, the corollary is proven.
\end{demostracion}

\begin{proposicion}{\label{0307046}} Let $n>1$ and
$y=\vg{y}{n}$ be a vector  of positive integers
with \tmcd\ equal to $1$.
There exists a  $n\times n$ matrix $Y$, whose entries are
non-negative integers, with determinant equal to $1$, one of
whose columns is $y$.
\end{proposicion}

\begin{demostracion}
Let $a\in\ZZ$ and denote by $E_{ij}(a)$, $i\neq j$, the
$n\times n$ elementary
matrix, which is equal to the identity matrix with its
$(i,j)$-entry replaced by $a$.
Let us write $y$ as a column vector.

Let us assume that all the
$y_j$ are multiple of one of them, say
$y_i$; then it must be $y_i=1$. Left
multiplications by
matrices $E_{ji}(-y_j)$ allow us to transform
the column vector $y$ into a column vector having all
entries equal to zero, except the $i$-th one which is equal
to $1$. Remark that all the elementary matrices we
have used have a {\em negative} entry out of the main
diagonal.

Let us assume that no entry of the column
vector $y$ divides all the others and
let
$y_i$ be the smallest of all these entries. By assumption,
there must be a $j$ such that, in the euclidean division,
$y_j=q_jy_i+r_j$ with $0<r_j<y_i$. Left multiplication of
$y$ by
$E_{ji}(-q_j)$ puts $r_j$ at the position $j$, leaving the
other entries unchanged. In this way, we get a new column
vector such that the \tmcd\ of its entries is $1$ and
the minimum of these
entries has strictly decreased.
Remark that, again, we have used an elementary matrix with a
negative integer entry out of the main diagonal.
If we repeat this process, after a finite
number of steps, we fall in the preceeding situation.

Summing-up: we have proven that, by left multiplication of
the
column vector $y$ by elementary  matrices having negative
integer entries out of the main diagonal, we arrive at a
matrix which is a column of the identity matrix. If we
denote by $Y$ the inverse of this product of elementary
matrices, we see that $Y$ is a product of elementary
matrices with positive integer entries out of the main
diagonal, so $Y$ is a matrix with non-negative integer
entries.
Let $I_{n}$ be the $n\times n$ identity matrix; then it is
clear that $Y=YI_{n}$ contains a column equal to $y$ and, of
course, $\det(Y)=1$. This proves the proposition.
\end{demostracion}

\begin{proposicion}{\label{0407041}} Let $Y$ be a
square matrix
with non-negative integer entries whose determinant is equal
to
$1$. Then $Y$ can be written as a product of
(by order): a finite
number of \mbu\ matrices, a permutation matrix and
another finite number of \mbu\ matrices.
\end{proposicion}

\begin{demostracion}
It is enough to show that $Y$ can be brought to a
permutation matrix by left and right multiplicaction
by \mbd\ matrices. If $q\in\ZZ_{+}$,
then $E_{ij}(-q)=E_{ij}(-1)^{q}$;
since $E_{ij}(-1)$ is a \mbd\ matrix, then $E_{ij}(-q)$ is a
product of \mbd\  matrices.

Let $y_{ij}$ be the smallest non-zero entry in $Y$. If all
the elements in the $j$-th column are multiple of $y_{ij}$,
we can get zeros in all the positions of this column,
except
$(i,j)$, by left multiplication by elementary matrices with
a negative entry out of the main diagonal.
In this case, the
fact that $\det(Y)=1$ implies $y_{ij}=1$. If some entry in
the $j$-th column is not a multiple of $y_{ij}$, say
$y_{lj}$, and if $y_{lj}=qy_{ij}+r_{lj}$ is the euclidean
division, then left multiplication by $E_{jl}(-q)$
puts $r_{lj}$ at the position
$(l,j)$ so the smallest non-zero entry of $Y$ has strictly
decreased. If we repeat this process, it is clear that we
must arrive to the first case after a finite number of
steps. The end of this process is a matrix $Y'_{1}$ with a
$1$ at one position (denote it again by $(i,j)$) and zeros
everywhere else in the $j$-th column. Moreover, $Y'_{1}$ is
the result of left multiplying $Y$ by \mbd\ matrices.

By symmetry, it is clear that we can get a new matrix
$Y_{1}$, obtained from $Y'_{1}$ by right multiplication by
\mbd\ matrices, and having $1$ at the $(i,j)$ position and
zeros everywhere else in the $i$-th row and the $j$-th
column. This is the basic argument of our proof.

We may repeat the argument for the
submatrix of $Y_{1}$ obtained by deleting the $i$-th row and
the $j$-th column, but seeing the operations in the whole
$Y_{1}$. This does not affect the form of $Y_{1}$. The very
end of the process is a matrix $Y_{p}$, which is a
permutation of the rows of the identity matrix, and which is
obtained from $Y$ by left and right multiplication by \mbd\
matrices. This proves the proposition.
\end{demostracion}

\begin{lema}{\label{0407042}} In the situation of theorem
\ref{0307044},
let $A$ be
the matrix whose row vectors are
$\{\mathbf{a}_{1},\ldots,\mathbf{a}_{m}\}$.
For all
$x=\vg{x}{m}\in\RR_0^m\setminus\{\mathbf{0}\}$ one has
$$
0\neq \sum_{i=1}^m  x_i\mathbf{a}_i\, .
$$
\end{lema}

\begin{demostracion}
Let $f\dps \RR^{n}\to\RR$ be  a linear function such that
$f(\mathbf{a}_i)>0$, $\forall i=1,\ldots,m$.
There always exists such a function: we give an example.
Let $g$ be an affine funcion such that $H$ has the equation
$g=0$ and $g(\mathbf{0})<0$. Then, for all $i=1,\ldots,m$
one has
$0=g(\mathbf{a}_{i})=g(\mathbf{0})+\overrightarrow{g}(\mathbf{a}_{i})$, so
$\overrightarrow{g}(\mathbf{a}_{i})>0$ and we may
thake $f=\overrightarrow{g}$.

Then,
$$
f\left(  \sum_{i=1}^m
x_i\mathbf{a}_i \right)=
\sum_{i=1}^m  x_i f(\mathbf{a}_i)>0\, ,
$$
which proves the lemma.
\end{demostracion}

\begin{nota}{\label{0407043}}
{\sc Proof of theorem \ref{0307044}:}

Let us assume that
$\mathbf{0}$ is the only point of $\Gamma(\Delta)$ belonging
to $-\RR_{0}^{n}$.
By lemma
\ref{0407042}, for every
$x\in\RR_0^m\setminus\{\mathbf{0}\}$, the vector $xA$ is
different from zero. By assumption, it cannot be
$xA\leq 0$. By corollary \ref{0307045}, there must exist
a vector $y\in\ZZ_+^n$ such that, written as a column
vector,
$Ay>0$. We may assume that the \tmcd\ of the entries of $y$
is $1$. By proposition \ref{0307046}, there exists a matrix
$Y$ with non-negative integer entries and determinant equal
to $1$ such that
$y$ is one of its columns. Therefore, one of the columns of
$AY$  is $Ay>0$, which implies by remark \ref{0407044} that
$AY$ can be brought to the \pc\ by a finite sequence of
\mbus. By proposition \ref{0407041}, we obtain $AY$ from $A$
by applying to $A$ a finite sequence of \mbus, then a
permutation of the columns and then another finite sequence
of \mbus. It is evident that the permutation of the columns
plays no role: if the matrix with the permuted columns can
be
brought to the \pc, also the original one. This proves the
theorem.

\end{nota}

\section{Applications} The main application of theorem
\ref{0307044} we consider here lies is to
the resolution of equations of the form: a
Weierstra\ss\
polynomial equal to zero. Let $k$ be an algebraically closed
field of
characteristic zero, $\mathbf{x}=\vg{x}{n}$ a collection of
indeterminates, $R=k\[[\mathbf{x}]]$ the corresponding ring
of power series, and let
$$
P(\mathbf{x},z)=z^{m}+h_{1}(\mathbf{x})z^{m-1}+\cdots+h_{m-1}(\mathbf{x})z+h_{m}(\mathbf{x})
\in k\[[\mathbf{x}]][z] \, , \quad n>1
$$
be an irreducible Weierstra\ss\ polynomial; the object to
study is the
equation $P(\mathbf{x},z)=0$. Let $D\in k\[[\mathbf{x}]]$ be
the discriminant of $P$ with respect to $z$; the
Jung-Abhyankar theorem (c.f. \cite{Abhyankar}) asserts that,
if $D$ is of the form $\mathbf{x}^{\mathbf{a}}U(\mathbf{x})$
with $U(\mathbf{x})\in k\[[\mathbf{x}]]$,
$U(\mathbf{0})\neq 0$, then the roots of $P(\mathbf{x},z)=0$
are a full set of conjugate Puiseux power series in the
variables $\mathbf{x}$. Here, $\mathbf{x}^{\mathbf{a}}$
means $\mathbf{x}^{\mathbf{a}}=x_{1}^{a_{1}}\cdots
x_{n}^{a_{n}}$, where
$\mathbf{a}=\vg{a}{n}\in\ZZ_{0}^{n}\setminus\{\mathbf{0}\}$
When the discriminant has this very special form, we say
that it is a {\em \ncd}\index{\ncd}.

In general, the roots are not Puiseux power series in
$\mathbf{x}$. However, we can say something very important
about them, namely

\begin{teorema}{\label{0507043}} The roots of
$P(\mathbf{x},z)=0$ are
 power series belonging to a ring
$k((x_{n}^{1/p}))\cdots((x_{i+1}^{1/p}))\[[x_{1}^{1/p},\ldots,x_{i}^{1/p}]]$
whose
\nds\ are contained in a pyramid $\Gamma(\Delta)$ such that
$\Gamma(\Delta)\cap(-\RR_{0}^{n})=\{\mathbf{0}\}$.
\end{teorema}

\vspace{0.3cm}

We will prove the theorem through several remarks.

\begin{nota}{\label{0607041}} We take the lexicographic
order in the sense (c.f. \cite{ZS2}, page 50):
if $\mathbf{a}\, , \mathbf{b}\in\RR^{n}$ then
$\mathbf{a}\mlex\mathbf{b}$ if and only if the first
component (from left to right) of $\mathbf{a}$ which is
different of the corresponding in $\mathbf{b}$
is strictly smaller.

Let us observe that a \mbu\ $\varphi_{ij}$ of the type
$\mathbf{a}\to\mathbf{a}M_{ij}$ with $i<j$ preserves the
lexicographic order, so it is an ordered automorphism
of $\RR^{n}$ (endowed with the lexicographic order).  In
fact,
$$
\varphi_{ij}(a_{1},\ldots,a_{i},\ldots,a_{j},\ldots,a_{n})=
(a_{1},\ldots,a_{i},\ldots,a_{i}+a_{j},\ldots,a_{n})\, ,
$$
so, if $\mathbf{a}\mlex\mathbf{b}$, then:
\begin{enumerate}

\item  If
$\mlex$ is decided before the
position $j$, it is evident that
$\varphi(\mathbf{a})\mlex\varphi(\mathbf{b})$.

\item
If $\mlex$
is decided at the position $j$, this means that
$a_{l}=b_{l}$, $\forall l$, $1\leq l\leq j-1$ and
$a_{j}<b_{j}$. Therefore, $a_{i}+a_{j}<b_{i}+b_{j}$, hence
$\varphi(\mathbf{a})\mlex\varphi(\mathbf{b})$.

\item
If $\mlex$
is decided at a position $l$ after $j$, this means that
all the components of
$\mathbf{a}$ until the
$(l-1)$-th coincide with the corresponding in
$\mathbf{b}$, so the same happens with
$\varphi(\mathbf{a})$ and $\varphi(\mathbf{b})$.
Since $a_{l}<b_{l}$
then
$\varphi(\mathbf{a})\mlex\varphi(\mathbf{b})$.

\end{enumerate}
We call this an {\em order-preserving \mbu}. Notice that
the corresponding \mbd\ $\varphi_{ij}^{-1}$ is also
order-preserving.
\end{nota}

\begin{notas}{\label{0507044}} Let
$\Lambda\subset\ZZ_{0}^{n}$ be a non-empty cloud of points;
we call the transform of $\Lambda$ by
a \mbu\ (or a \mbd) the set of the transforms of all the
points of
$\Lambda$.

\parnot{\label{0507044a1}}
Any \mbu\  $\varphi$ keeps
$\ZZ_{0}^{n}$, that is,
$\varphi(\ZZ_{0}^{n})\subset\ZZ_{0}^{n}$.  Therefore, if
$\mathbf{a}\,
,
\mathbf{b}\in\ZZ_{0}^{n}$ are two points such that
$\mathbf{b}\in\mathbf{a}+\ZZ_{0}^{n}$
then
$\varphi(\mathbf{b})\in\varphi(\mathbf{a})+\ZZ_{0}^{n}$.

\parnot{\label{0507044a11}}
Let us assume that $\varphi=\varphi_{pq}$,
$p<q$, is an order-preserving \mbu.
Let $\mathbf{a}\, , \mathbf{b}\in\ZZ_{0}^{n}$ and
write $\mathbf{a}'=\varphi(\mathbf{a})$,
$\mathbf{b}'=\varphi(\mathbf{b})$. Let us assume that
there exists
an index $j$, $1\leq j\leq n$ such that
$b_{i}\geq a_{i}$, $\forall i=1,\ldots,j$; then
$b'_{i}\geq a'_{i}$, $\forall i=1,\ldots,j$. In fact, the
result is clear by \ref{0507044}.\ref{0507044a1} if $j=n$, so let us
assume that $j<n$. Let
$\mathbf{c}=(b_{1}-a_{1},\ldots,b_{j}-a_{j},0,\ldots,0)$
and
$\mathbf{d}=(0,\ldots,0,b_{j+1}-a_{j+1},\ldots,b_{n}-a_{n})$;
then $\mathbf{b}=\mathbf{a}+\mathbf{c}+\mathbf{d}$.
Since $\mathbf{a}+\mathbf{c}\in\mathbf{a}+\ZZ_{0}^{n}$, then
$\varphi(\mathbf{a}+\mathbf{c})\in\varphi(\mathbf{a})+\ZZ_{0}^{n}$.
On the other hand, since $p<q$, the first
$j$ components
of $\varphi(\mathbf{d})$ are zero, so the conclusion
is clear. It is obvious that the same happens if we replace
$\varphi_{pq}$  by the composition of a finite sequence of
order-preserving \mbus.

\parnot{\label{0507044a12}}
For any $j=1,\ldots,n$ and any $\mathbf{v}\in\ZZ_{0}^{n}$,
we write
$\Lambda_{j}(\mathbf{v})=\{\mathbf{v}'\in\ZZ_{0}^{n}\mid
v'_{i}\geq v_{i}\, , \forall i=1,\ldots,j\}$;
 then
$\Lambda_{1}(\mathbf{v})\supset\cdots\supset\Lambda_{n}(\mathbf{v})$.
By \ref{0507044}.\ref{0507044a11}, for every
order-preserving \mbu\ $\varphi$, every
$\mathbf{v}\in\ZZ_{0}^{n}$ and every
$j=1,\ldots,n$, one has
$\varphi(\Lambda_{j}(\mathbf{v}))\subset\Lambda_{j}(\varphi(\mathbf{v}))$.
Moreover, if $\Phi$ is a composition of  a finite number of
order-preserving \mbus, then
$\Phi(\Lambda_{j}(\mathbf{v}))\subset\Lambda_{j}(\Phi(\mathbf{v}))$.

\parnot{\label{0507044a2}} Since
$\emptyset\neq\Lambda\subset\ZZ_{0}^{n}$, there is
a minimum-lex $\mathbf{u}\in\Lambda$, so
$\Lambda\subset\Lambda_{1}(\mathbf{u})$. We will now prove
that there
exists a finite sequence of order-pre\-ser\-ving \mbus\ such
that, calling
$\Phi$ the composition of all of them, one has
$\Phi(\Lambda)\subset\Phi(\mathbf{u})+\ZZ_{0}^{n}$.
If $\Lambda\subset\Lambda_{n}(\mathbf{u})$ there is nothing
to prove, so we assume this is not the case.
Let $j$ be the smallest index such that
$\Lambda\not\subset\Lambda_{j}(\mathbf{u})$; then,
necessarily $j>1$. By the minimality of $j$, for every
$\mathbf{u}'\in\Lambda$
one must have $u'_{i}\geq u_{i}$, $\forall i=1,\ldots,j-1$
Moreover, if
$\mathbf{u}'\in\Lambda\setminus\Lambda_{j}(\mathbf{u})$ then
$u'_{j}<u_{j}$, so $u_{j}>0$.
Since $\mathbf{u}\in\Lambda$ is the minimum-lex, for every
$\mathbf{u}'\in\Lambda\setminus\Lambda_{j}(\mathbf{u})$
there must exist an
index
$i<j$
such that $u'_{i}>u_{i}$. Let $i_{1}<j$ be the smallest
index
such that there exists $\mathbf{u}'\in\Lambda\setminus\Lambda_{j}(\mathbf{u})$
satisfying $u'_{i_{1}}>u_{i_{1}}$.
Let $\Phi_{1}$ be the
composition of $u_{j}$ \mbus\ equal to $\varphi_{i_{1}j}$;
then
$\Phi_{1}(\mathbf{u})=(u_{1},\ldots,u_{i_{1}},\ldots,u_{j}u_{i_{1}}+u_{j},\ldots,u_{n})$
and, for every $\mathbf{u}'\in\Lambda\setminus\Lambda_{j}(\mathbf{u})$
with $u'_{i}>u_{i}$,
$\Phi_{1}(\mathbf{u}')=(u'_{1},\ldots,u'_{i_{1}},\ldots,u_{j}u'_{i_{1}}+u'_{j},\ldots,u'_{n})$.
Since $u'_{i_{1}}>u_{i_{1}}$, then
$u_{j}(u'_{i_{1}}-u_{i_{1}})\geq u_{j}\geq u_{j}-u'_{j}$,
hence $u_{j}u'_{i_{1}}+u'_{j}\geq
u_{j}u_{i_{1}}+u_{j}$, so
$\Phi_{1}(\mathbf{u}')\in\Lambda_{j}(\Phi_{1}(\mathbf{u}))$.

Now, $\Phi_{1}(\mathbf{u})$ is the minimum-lex of
$\Phi_{1}(\Lambda)$ and, by
\ref{0507044}.\ref{0507044a12},
$\Phi_{1}(\Lambda_{i}(\mathbf{u}))=\Lambda_{i}(\Phi_{1}(\mathbf{u}))$, $\forall i=1,\ldots,n$. By
\ref{0507044}.\ref{0507044a11},
$\Phi_{1}(\Lambda)\subset\Lambda_{i}(\Phi_{1}(\mathbf{u}))$,
$\forall i=1,\ldots,j-1$. One must not forget that
$\Phi_{1}$ leaves invariant the first $j-1$ components of
every vector. Therefore, if
$\Phi_{1}(\Lambda)\not\subset\Lambda_{j}(\Phi_{1}(\mathbf{u}))$,
there exists a smallest index $i_{2}$
such that there exists $\mathbf{u}'\in\Phi_{1}(\Lambda)\setminus\Lambda_{j}(\Phi_{1}(\mathbf{u}))$
satisfying $u'_{i_{2}}>u_{i_{2}}$. Necessarily $i_{2}>i_{1}$
and we proceed as before, and so on. It is then clear that
there exists a finite sequence of order-preserving \mbus,
whose composition $\overline{\Phi}$ is such that
$\overline{\Phi}(\Lambda)\subset\Lambda_{i}(\overline{\Phi}(\mathbf{u}))$,
$\forall i=1,\ldots,j$. If there is still a $j_{1}$ such
that
$\overline{\Phi}(\Lambda)\not\subset\Lambda_{j_{1}}(\overline{\Phi}(\mathbf{u}))$,
then $j_{1}>j$ and we proceed as before, and so on. This
proves our assertion.
\end{notas}

\begin{nota}{\label{0507045}}{\sc Proof of theorem
\ref{0507043}}. Let $\Lambda$ be the \nd\ of the
discriminant $D$ of $P(\mathbf{x},z)$; by
\ref{0507044}.\ref{0507044a2} there exists a finite sequence
of order-preserving \mbus\ such that, calling $\Phi$ their
composition,
 $\Phi(\Lambda)\subset \mathbf{a}+\ZZ_{0}^{n}$ where
$\mathbf{a}\in\Phi(\Lambda)$. We make these \mbus\ to act
upon $P(\mathbf{x},z)$ and denote by $Q(\mathbf{x},z)$
the transform of $P(\mathbf{x},z)$ by $\Phi$. The
discriminant $D'$ of $Q(\mathbf{x},z)$ is just the transform
of $D$ because $D$ is a polynomial in the coefficients of
the equation. Moreover, $D'$ is a \ncd, hence the roots of
$Q=0$ are all ordinary Puiseux power series, say
with common denominator $p$ of the exponents, because
every
irreducible factor of $Q(\mathbf{x},z)$ has a discriminant
which is a \ncd.  If we come back to the beginning by
applying the corresponding sequence of \mbds, the region
containing the \nd\ of the roots of $Q=0$, namely the \pc,
obviously goes to a pyramid $\Gamma(\Delta)$ such that
$\Gamma(\Delta)\cap(-\RR_{0}^{n})=\{\mathbf{0}\}$.
Since all the \mbus\ are of the form $\varphi_{lj}$  with
$l<j$, we denote by $i$ the minimum of all the indices $l$
of these \mbus, then $\Phi$ leaves invariant the first $i$
coordinates of every point, so the same happens with
$\Phi^{-1}$. Therefore,
every
monomial
$x_{1}^{a_{1}/p}\cdots x_{i}^{a_{i}/p}x_{i+1}^{a_{i+1}/p}\cdots x_{n}^{a_{n}/p}$
occurring in a root of $Q$
evolves in a way such that the
exponents
$a_{1}/p,\ldots,a_{i}/p$ remain unchanged. Therefore, if we
fix a root $\varrho$ of $Q=0$, fix
$a_{1}/p,\ldots,a_{i}/p$ and write
$\varrho'=x_{1}^{a_{1}/p}\cdots x_{i}^{a_{i}/p}
\varrho''(x_{i+1}^{a_{i+1}/p},\ldots,x_{n}^{a_{n}/p})$,
with
$\varrho''(x_{i+1}^{a_{i+1}/p},\ldots,x_{n}^{a_{n}/p})\in
k\[[x_{i+1}^{1/p},\ldots,x_{n}^{1/p}]]$, for the sum of
all the terms  of the root whose monomials start by $x_{1}^{a_{1}/p}\cdots x_{i}^{a_{i}/p}$, the
transform
of $\varrho'$ by $\Phi^{-1}$ produces a power
series
$x_{1}^{a_{1}/p}\cdots
x_{i}^{a_{i}/p}
\varrho''_{1}(x_{i+1}^{a_{i+1}/p},\ldots,x_{n}^{a_{n}/p})$
where $\varrho''_{1}$ has possibly negative exponents.
Since $(1/p)\cdot\ZZ_{0}^{n}$ is lexicographically well-ordered,
so it is $\Phi^{-1}((1/p)\cdot\ZZ_{0}^{n})$, hence the
transform of $\varrho$ by $\Phi^{-1}$ belongs to
$k((x_{n}^{1/p}))\cdots((x_{i+1}^{1/p}))\[[x_{1}^{1/p},\ldots,x_{i}^{1/p}]]$,
which proves the theorem.

\end{nota}

\vspace{0.3cm}

When $n=2$ there is much more to say, namely:

\begin{nota}{\label{0607042}} In our joint
paper (cf. \cite{SV1}),
we prove the following for $n=2$:
\begin{enumerate}

\item The theorem \ref{0507043} {\em without using the
Jung-Abhyankar theorem}.

\item The Jung-Abhyankar theorem from the fact that the
\nds\ of the roots lie in a pyramid
$\Gamma(\Delta)$ such that
$\Gamma(\Delta)\cap(-\RR_{0}^{2})=\{\mathbf{0}\}$.

\end{enumerate}
This shows that the Jung-Abhyankar theorem in dimension
$2$ can
be proven by linear algebra techniques, without having
resource to more sophisticated algebraic material.
Moreover, in this case, the Jung-Abhyankar theorem is
equivalent to the fact that the roots of the equation lie in
a pyramid satifying the conditions of  theorem
\ref{0307044}.

\end{nota}

\section{Short remarks on computations} The
explicit computations are a consequence, more or
less obvious, of the convex calculus and the optimization
of a linear function on a polyhedron by the simplex method.

The point of departure will be always the list of points
$\{\mathbf{a}_{1},\ldots,\mathbf{a}_{m}\}$, all different
from
$\mathbf{0}$, generating, either a hyperplane not passing
through the origin, or $X=\RR^{n}$. We add
$\mathbf{0}$ to the list, and write
$E=\{\mathbf{0},\mathbf{a}_{1},\ldots,\mathbf{a}_{m}\}$; in
both cases $E$ generates the whole affine space. We denote
by
$A$ the matrix whose row vectors are
$\{\mathbf{a}_{1},\ldots,\mathbf{a}_{m}\}$.

\begin{nota}{\label{0407046}} By elementary linear calculus
(c.f. Vicente, J.L., loc. cit.), the $(n-1)$-dimensional
faces of the polyhedron $[E]$ are produced by the
following algorithm: we pick all the subsets
of $E$ consisting of $n$ \ai\ points, and find the
hyperplane determined by them; then we drop repetitions and
keep only those hyperplanes leaving all the points of $E$ in
an only \se. This algorithm is not the best possible,
but improvements are out of the scope of this paper.
If $C=\{H_{1},\ldots,H_{p}\}$ is the list of
faces, then we get the vertices of $[E]$ by the following
algorithm: we pick all the subsets of $C$
consisting of $n$ hyperplanes whose intersection is an only
point, find the point,  drop repetitions and the
remaining ones are the
vertices. It is clear that $E$ defines a pyramid $\Gamma$ if
and only if $\mathbf{0}$ is a vertex of $[E]$. In this case,
the faces of $\Gamma$ are those $H_{i}$ passing through
$\mathbf{0}$.
For instance,
if we start from the points
$$
\mathbf{a}_{1}=(4,7,-9)\, , \;
\mathbf{a}_{2}=(5,7,-8)\, , \;
\mathbf{a}_{3}=(3,5,-9)\, , \;
\mathbf{a}_{4}=(4,0,-1) \, ,
$$
which generate $\RR_{0}^{n}$, the faces are
$$
\begin{array}{ll}
7\,x_{{1}}-13\,x_{{2}}-7\,x_{{3}} = 0 &
18\,x_{{1}}-9\,x_{{2}}+x_{{3}}=0 \\
-7\, x_{{1}}-27\,x_{{2}}-28\,x_{{3}}= 0 &
5\,x_{{1}}+33\,x_{{2}}+20\,x_{{3}} = 0 \\
19-2\,x_{{1}}+x_{{2}}+2\,x_{{3}} =0 &
32-7\,x_{{1}}+5\,x_{{2}}+4\,x_{{3}}=0
\end{array}
$$
and the vertices are
$E=\{\mathbf{0},\mathbf{a}_{1},\mathbf{a}_{2},\mathbf{a}_{3},\mathbf{a}_{4}\}$;
therefore
the points define a pyramid. The faces are normalized in the
sense that all the points of $E$ make their linear
equations $\geq 0$. The faces of the pyramid are the first
four and the edges are the positive \tsrs\ determined by the
four given points.
\end{nota}

\begin{nota}{\label{0507041}} It is not difficult to know
whether the pyramid $\Gamma$ satisfies the condition
$\Gamma\cap(-\RR_{0}^{n})=\{\mathbf{0}\}$ or not.
Let
$\Lambda=\vg{\lambda}{m}$ be a row of variables and let
$m_{i}$ be the element in the $i$-th column of the matrix
$\Lambda A$; then
$\Gamma\cap(-\RR_{0}^{n})\neq\{\mathbf{0}\}$ if and only if
there is
a feasible solution to the set of linear constraints
$$
1+\sum_{i=1}^{m}m_{i}=0\, ,\quad m_{i}\leq 0\, , \quad
\lambda_{i}\geq 0\, , \quad i=1,\ldots,m\, .
$$
The existence of a feasible solution can obviously be
decided by the simplex method. In the preceeding example,
the feasible solution does not exist, so $\Gamma\cap(-\RR_{0}^{n})=\{\mathbf{0}\}$.
\end{nota}

\begin{nota}{\label{0507042}} It is also easy to
find a positive solution of the system of inequalities
$Ay>0$, where $y$ is a column of variables. If $m_{i}$ is
the element in the $i$-th row of $Ay$, we can easily get a
positive
solution of $Ay>0$ by minimizing any of the coordinate
functions on the set of constraints
$m_{i}\geq 1\, , y_{i}\geq 1$, $i=1,\ldots,m$. In our
example, minimizing $y_{1}$ by the simplex method will
produce the point $(1,7/5,1)$, so the integer solution
$(5,7,5)$.
\end{nota}

\vspace{0.3cm}

The remaining computations
to bring $A$ to the \pc\ are straightforward matrix
operations.

\bibliographystyle{plain}
\bibliography{kk}

\end{document}